\documentclass[10pt,reqno]{amsart}
\usepackage{amsmath,amsthm,amssymb,amscd,url,enumerate}
\usepackage[dvipsnames]{xcolor}
\usepackage{bm} 
\usepackage{graphicx}
\usepackage{afterpage}
\usepackage{tikz}
\usepackage[all]{xy}
\usepackage{mathrsfs}
\usepackage[margin=1in]{geometry}
\usepackage[ruled,vlined,linesnumbered,algosection]{algorithm2e}
\usepackage{algpseudocode}
\usepackage[inline]{enumitem}
\usepackage{aliascnt} 
\usepackage[colorlinks=true,linkcolor=blue,citecolor=RedViolet,urlcolor=RedViolet]{hyperref}

\newcommand{\TITLE}{Factoring using multiplicative relations modulo $n$: a subexponential algorithm inspired by the index calculus}
\newcommand{\TITLERUNNING}{\TITLE}
\newcommand{\DATE}{\today}

\newaliascnt{theorem}{algocf}
\theoremstyle{plain} 
\newtheorem{theorem}[theorem]{Theorem}

\newtheorem{hypothesis}[theorem]{Hypothesis}

\theoremstyle{definition}

\theoremstyle{remark}
\newtheorem{remark}[theorem]{Remark}

\numberwithin{theorem}{section}

  {\end{list}}

  {\end{list}}

\newcommand{\tightoverset}[2]{%
  \mathop{#2}\limits^{\vbox to -.5ex{\kern-1.05ex\hbox{$#1$}\vss}}}

\newcommand{\NN}{\mathbb{N}}

\newcommand{\QQ}{\mathbb{Q}}

\newcommand{\ZZ}{\mathbb{Z}}

\newcommand{\GCD}{{\operatorname{GCD}}}
\renewcommand{\gcd}{{\operatorname{gcd}}}

\newcommand{\MOD}[1]{~(\textup{mod}~#1)}
\renewcommand{\pmod}{\MOD}

\newcommand{\ord}{\operatorname{ord}}

\newcommand{\poly}{\operatorname{poly}}

\newif\ifcomments
\commentstrue  

\definecolor{myyellow}{rgb}{1.0, 0.75, 0.0}
\definecolor{mygreen}{rgb}{0.35, 0.71, 0.0}

\ifcomments
\newcommand{\KS}[1]{\textcolor{blue}{{\sf (ToDo:} {\sl{#1})}}}
\else
\newcommand{\KS}[1]{}
\fi

\title[\TITLERUNNING]{\vspace*{-1.3cm} \TITLE}

\author{Katherine E. Stange}

\date{\DATE}
\thanks{Katherine E. Stange was supported by NSF-CAREER CNS-1652238.}
\address{%
Department of Mathematics, University of Colorado,
Campux Box 395, Boulder, Colorado 80309-0395}
\email{kstange@math.colorado.edu}

\begin{document}

\maketitle

\begin{abstract}
	We demonstrate that a modification of the classical index calculus algorithm can be used to factor integers.  More generally, we reduce the factoring problem to finding an overdetermined system of multiplicative relations in any factor base modulo $n$, where $n$ is the integer whose factorization is sought.  The algorithm has subexponential runtime $\exp(O(\sqrt{\log n \log \log n}))$ (or $\exp(O( (\log n)^{1/3} (\log \log n)^{2/3} ))$ with the addition of a number field sieve), but requires a rational linear algebra phase, which is more intensive than the linear algebra phase of the classical index calculus algorithm.  The algorithm is certainly slower than the best known factoring algorithms, but is perhaps somewhat notable for its simplicity and its similarity to the index calculus.
\end{abstract}

\section{Introduction}

The index calculus (whose roots date to Kraitchik \cite{Kraitchik} and which was rediscovered and developed more recently by Adleman \cite{Adleman} and Western and Miller \cite{WesternMiller}, among others) computes the discrete logarithm of a residue $h$ modulo $p$ with respect to a base $g$ by collecting relations of the form
\[
	g^{x} = \prod_{i=1}^b p_i^{e_{i}} \pmod{p},
\]
in terms of a factor base of primes $p_i$.  Choosing subexponentially many primes, it takes subexponential time to find as many relations as there are primes.  This leads to a linear algebra system, whose solution gives the discrete logarithms of the primes, which can then be used to find the discrete logarithm of $h$.

It has long been observed that factoring algorithms and discrete logarithm algorithms have many similarities.  In the present note, we demonstrate a factoring algorithm that is especially closely inspired by the index calculus algorithm.  The fundamental idea is easy to state:  we attempt to use the index calculus algorithm modulo an integer $n$ with unknown factorization (instead of a prime $p$).  That is, with respect to a factor base of primes, collect relations of the form
\[
	g^{x} = \prod_{i=1}^b p_i^{e_{i}} \pmod{n}.
\]
In this case, the exponent relations live modulo the Euler totient $\varphi(n)$, or, more precisely, modulo the order of the chosen base $g$.  We do not know this order, and in fact its discovery would typically lead to a factorisation of $n$.  Although the discrete logarithms of a factor base need not exist in this situation, multiple relations will tend to lead to a linear system that is overdetermined and has no solution, unless working modulo the order of $g$.  This leads to a linear algebra method to extract the order of $g$. 

It turns out that more generally, any method of finding an overdetermined system of multiplicative relations between elements of a factor base modulo $n$ will lead to a method of factorization.  We give a simple heuristic reduction from factoring to finding such a system of multiplicative relations:

\begin{theorem}[Introductory form of Theorem~\ref{thm:reduc}]
	\label{thm:reduc-intro}
	Let $n$ be an integer, and let $b$ be a function of $n$. 
	Let $\mathcal{B}$ be a factor base of $b$ residues modulo $n$. 
	Suppose an oracle provides random multiplicative relations of size $\le O(\log n)$ amongst the residues $\mathcal{B}$.  
	Then there is a Monte Carlo algorithm to determine the multiplicative order of residues modulo $n$ (and hence factor $n$), whose runtime is polynomial in $b$ and $\log n$, and which requires at most $O(b)$ calls to the oracle.  
	Furthermore, under Hypothesis~\ref{heur:sub} (Section~\ref{sec:reduc}), the probability of success approaches $1 - 1/\zeta(c+1)$ where a total of $b+c$ calls are made to the oracle.
\end{theorem}

The notation $\zeta$ denotes the Riemann zeta function.
Although the algorithm may fail (it may return a nontrivial multiple of the order of $g$, or return $0$), a simple and informal implementation using SageMath \cite{SageMath} indicates the algorithm succeeds with very high probability.  Hypothesis~\ref{heur:sub} concerns the probability that a random selection of elements of a lattice will generate that lattice, or, more precisely, that the same is true for a restriction to a one-dimensional subspace.  The probability $1 - 1/\zeta(k)$ is more familiar as the probability that $k$ random integers have no non-trivial common factor.

The runtime of the new index-calculus-like algorithm is slower than that of the textbook index calculus, because the linear algebra step involves arithmetic over $\QQ$.  Its runtime is $\exp(O( (\log n)^{1/2} (\log \log n)^{1/2} ))$.  
The textbook index calculus described above has been improved by the addition of the number field sieve in the precomputation (relation-finding) phase \cite{Gordon}.  The methods of \cite[Section 3.1]{Gordon} can be adapated to find relations modulo $n$, which, when combined with Theorem~\ref{thm:reduc}, leads to a version of the present algorithm which runs in time $\exp(O( (\log n)^{1/3} (\log \log n)^{2/3} ))$.  We do not pursue this in greater detail here.

The literature for factoring includes a panoply of approaches.  As far as index-calculus-style approaches to factoring, another family of approaches uses an index calculus within the class group of a quadratic field \cite{Lenstra-fast, Lenstra-Pomerance, Seysen}.  The analysis of Section~\ref{sec:reduc} has some commonalities with the Hafner-McCurley rigorous class group algorithm \cite{HafnerMcCurley}, and their methods may be applicable to making the present algorithm rigorous.  Although the present algorithm is closely related to the index calculus and other relation-finding approaches such as the quadratic and number fields sieves, the author has been unable to find this particular variation in the literature.
For background on the index calculus, the number field sieve and factoring algorithms in general, see \cite{CrandallPomeranceBook}.  

\subsection*{Accompanying code}

A toy implementation is available at 
\begin{center}\url{https://github.com/katestange/index-factor}
\end{center}

\subsection*{Acknowledgements}
Thank you to Carl Pomerance and Sam Wagstaff for sharing their expertise in the factoring literature, and to Sam Wagstaff and Renate Scheidler for helpful comments on an earlier draft.  Thank you to Karthekeyan Chandrasekaran for help with a reference.  And finally, the author is grateful to the students in her Fall 2022 Introduction to Coding Theory and Cryptography at the University of Colorado Boulder for their enthusiasm and engagement, which led, indirectly, to this paper.

\section{The algorithm}

The algorithm is actually one for finding the multiplicative order of $g$ modulo $n$.  It is well-known that this leads to a factorization of $n$.  For example, Shor's quantum factoring algorithm determines this order $r$ for random values of $g$ until it finds an order $r$ which is even and for which $g^{r/2} \neq -1 \pmod{n}$ \cite{Shor}.  In this case, $\gcd(g^{r/2} \pm 1, n)$ should be a non-trivial factor.   More generally, including odd orders, see \cite{Ekera}.

By the \emph{period problem} for an integer $n$, we shall mean the problem of computing the multiplicative order of a given residue modulo $n$.  The \emph{factorization problem} for an integer $n$ is the problem of giving its unique prime factorization over the integers.  It is known that these are equivalent.

\begin{theorem}[{\cite[Theorem 4]{Miller}}]
	Under the Extended Riemann Hypothesis, the period problem and the factorization problem are polynomial-time equivalent.
\end{theorem}

The Extended Riemann Hypothesis is used to guarantee the existence of a small residue which generates the multiplicative group; in practice, choosing a residue at random is likely to result in a generator very quickly.

We now describe the algorithm.  Let $n$ be an odd positive integer, and let $g$ be a residue modulo $n$.  Let $\mathcal{B}$ be a factor base consisting of the primes $p_1, \ldots, p_b$ less than or equal to a bound $B$, chosen depending upon $n$ (so $b = \#\mathcal{B}$).  

For random integers $0 < x_j < n$, compute $g^{x_j} \pmod{n}$ and attempt to factor the resulting smallest positive residue in the factor base.  For each relation obtained in this way, being of the form
\[
	g^{x_j} = \prod_{i=1}^b p_i^{f_{j,i}} \pmod{n},
\]
store the vector $\mathbf{f}_j = (f_{j,i})_{i=1}^b \in \ZZ^b$ and the value $x_j$.
Collect these vectors until we have more than $b$ vectors; say $\mathbf{f}_1, \ldots, \mathbf{f}_{b+c}$.  (Taking $c=10$ or so should generally suffice, but this is a parameter of the algorithm.)

We should now find integer relations between the vectors $\mathbf{f}_j$.  These relations can be found as elements of the basis of the right kernel of the $b \times (b+c)$ matrix whose columns are the relations, working over $\QQ$.  Using gcd computations, scale $c$ of these basis elements to have integral entries with no common factor, and call them $\mathbf{b}_1, \ldots, \mathbf{b}_c \in \ZZ^{b+c}$.  In other words,
\[
	\sum_{j=1}^{b+c} (\mathbf{b}_t)_j \mathbf{f}_j = \mathbf{0},
\]
for $t=1,\ldots,c$.
Collect the values
\[
	\alpha_t := \sum_{j=1}^{b+c} (\mathbf{b}_t)_j x_j.
\]
We return the $\gcd$ of this set of values.  We show below that this $\gcd$ is is with high probability the multiplicative order of $g$ modulo $n$, and otherwise an integer multiple of this multiplicative order.  The algorithm is given more formally in Algorithm~\ref{alg:factor}.

\begin{algorithm}
	\caption{\,Computing the multiplicative order of $g$ modulo $n$.} \label{alg:factor}
    \vspace{.2ex}
	\SetKwInOut{Input}{Input}
	\SetKwInOut{Output}{Output}
	\SetKwInput{Init}{Inititialization phase}
	\SetKwInput{Rel}{Relation finding phase}
	\SetKwInput{LA}{Linear algebra phase}
	\SetKwInput{GCD}{GCD computation phase}
    \Input{%
	    A positive integer $n$, and a positive integer $g < n$.
           }
    \Output{%
	    The multiplicative order of $g$ modulo $n$.
        }
    \BlankLine

    \Init{}

    \BlankLine

    Select a suitable $B \in \NN$, and let $\mathcal{B} = \{ p_1, p_2, \ldots, p_b \}$ be the \emph{factor base} consisting of all primes less than $B$.

    Select a suitable $c \in \NN$ (typically $c \sim 10$ is fine).

    \BlankLine
    \Rel{}

    \BlankLine
    Let $j = 1$.

    \While{ $j < b+c$ }{
	
	    Choose an integer $x$ randomly in the range $[1,\ldots, n]$.  (If $x$ has been drawn previously, draw again.)

	    Compute the smallest positive residue of $g^x$ modulo $n$.

	    Attempt to factor the residue of $g^x$ in terms of the factor base.

	    \If{ $g^x$ is factored, say $g^x = \prod_{i=1}^b p_i^{f_i} \pmod n$ }{

		    Set $\mathbf{f}_j = (f_1, \ldots, f_b)$.
		    Set $x_j = x$.
		    Increment $j$.
	    }
    }

    \BlankLine
    \LA

    \BlankLine
    Form the $b \times (b+c)$ integer matrix whose columns are the $\mathbf{f}_j$.

    Compute independent vectors $\mathbf{b}_1, \ldots, \mathbf{b}_c \in \QQ^{b+c}$ in the right kernel of this matrix (i.e., column combinations that vanish).

    \BlankLine
    \GCD

    \BlankLine
    Scale each basis element so that the entries are integers with no factor common to all entries.

    \For{ $t$ from $1$ to $c$ }{

	    Let $\alpha_t = \sum_{j=1}^{b+c} (\mathbf{b}_t)_j x_j$. \SetKwComment{comment}{/*}{*/}

    }

    Let $G$ be the greatest common divisor of all the $\alpha_1, \ldots, \alpha_c$.

   \Return{ $G$ }
\end{algorithm}

\subsection{Correctness}

We now briefly demonstrate the claim that the result is a multiple of $\ord(g)$, the multiplicative order of $g$.  For a more complete analysis of runtime and success probability, see the next section.
When one has a relation 
\begin{equation*}
	\sum_{j=1}^{b+c} (\mathbf{b}_t)_j \mathbf{f}_j = \mathbf{0},
\end{equation*}
this implies that 
\[
	\prod_{j=1}^{b+c} (g^{x_j})^{(\mathbf{b}_t)_j} = 1 \pmod{n}.
\]
In particular, this implies that
\[
	\alpha_t = \sum_{j=1}^{b+c} (\mathbf{b}_t)_j x_j = 0 \pmod{\ord(g)}.
\]
Hence, the multiplicative order $g$ modulo $n$ must divide every $\alpha_t$.  Therefore the final gcd of the algorithm is a multiple of $\ord(g)$.

\begin{remark}
Particularly when $g$ has small order, it is likely that the individual discrete logarithms $L_g(p_i)$ don't all exist.  However, the algorithm doesn't mind; the factorization is just a way to create multiplicative relations between values of $g^x$.  If $L_g(p_i)$ and $L_g(p_j)$ don't exist, but $L_g(p_i p_j)$ does, then $L_g(p_i)$ will only appear with the same coefficient as $L_g(p_j)$ in the linear algebra, effectively reducing the number of variables and increasing the size of the kernel.  But the correctness of the algorithm is not affected.
\end{remark}

\section{Factoring by relations modulo $n$}
\label{sec:reduc}

In order to analyse the algorithm's correctness and runtime, we will isolate the novel phase encompassing the linear algebra and gcd computations.  In fact, this phase of the algorithm above generalizes to give a heuristic reduction from factoring to finding multiplicative relations modulo $n$, which we will state in Theorem~\ref{thm:reduc} below.  This theorem will then imply that Algorithm~\ref{alg:factor} is correct and has a high probability of success, under Hypothesis~\ref{heur:sub} below.  This hypothesis is known to hold in some cases, such as when the size of the factor base is much smaller than $n$, so that in those cases there is a reduction from factoring $n$ to finding relations in such a factor base (under ERH).

If we have a factor base $\mathcal{B}$ of $b$ residues $a_1, \ldots, a_b$ modulo $n$, then by a \emph{multiplicative relation}, we shall mean a relationship
\[
	\prod_{i=1}^b a_i^{e_i} = 1
\]
for $e_i \in \ZZ$.  We give a name to the lattice of exponent vectors for all valid relations:
\[
	\Lambda_\mathcal{B} = 
	\left\{ \mathbf{e} = (e_i)_{i=1}^b : 
	\prod_{i=1}^b a_i^{e_i} = 1 
\right\}.
\]
This is a sublattice of $\ZZ^b$.  If the residues generate $(\ZZ/n\ZZ)^*$, then it will have covolume equal to $\varphi(n)$.  

Let $S_i$ be the $i$-th coordinate subspace, i.e. vectors whose entries are non-zero only in the $i$-th position.  The restriction $\Lambda_\mathcal{B}|_{S_i}$ of $\Lambda_\mathcal{B}$ to $S_i$ has covolume equal to the multiplicative order of $a_i$, denoted $\ord(a_i)$.  In particular, it is generated by a vector having $i$-th entry $\ord(a_i)$.  Equivalently, any generating set for $\Lambda_\mathcal{B}|_{S_i}$ will have entries with $\gcd$ equal to $\ord(a_i)$.

The \emph{size} of a relation will be equal to the logarithm of the $1$-norm $|\mathbf{e}_i|_1$ of its exponent vector.  (Thus relation vectors whose entries are $< n$ have size $O(\log n)$.)

\begin{hypothesis}[Hypothesis on sublattice generation]
	\label{heur:sub}
	Let $n > 0$ tend to $\infty$, and $b = f(n)$ be a given function of $n$.
	Let $\Lambda_\mathcal{B}' \subseteq \Lambda_\mathcal{B}$ be a lattice generated by $b+c$ relations randomly chosen from amongst those in $\Lambda_\mathcal{B}$ of size $O(\log n)$.  Then, asymptotically, the probability that $\Lambda_\mathcal{B}'|_{S_i} = \Lambda_\mathcal{B}|_{S_i}$ is equal to the probability that $c+1$ random integers (in the sense of natural density) share no common factor, i.e. $1 - 1/\zeta(c+1)$ where $\zeta$ is the Riemann zeta function.
\end{hypothesis}

If $b$ and $n$ satisfy the relationship $n \ge 8 b^{b/2}$ as $n$ tends to infinity, then taking $c = b+1$, it is known that the probability has a lower bound \cite[Theorem 1.1]{FonteinWocjan}.  By contrast, to apply this result to Algorithm~\ref{alg:factor}, we need the Hypothesis to hold when $b$ is subexponential in $\log n$.

To provide some evidence for the Hypothesis, we give a heuristic argument.
The index $[\Lambda_\mathcal{B}|_{S_i} : \Lambda_\mathcal{B}'|_{S_i}]$ is computed (for example in Algorithm~\ref{alg:factor}) as follows:  let $K$ be the right kernel of the matrix $M$ whose columns are the $b+c$ relations.  It has dimension $\ge c$, and with high probability the dimension is exactly $c$ (if the relations generate a sublattice of $\Lambda_\mathcal{B}$ of full rank; see \cite{FonteinWocjan} for a result that this is the case with constant probability once the relation size is large enough).  Let $M_i$ be the matrix obtained by deleting the $i$-th row of $M$.  Then its right kernel $K_i \supseteq K$ (of dimension $c+1$ with high probability) gives linear combinations of the original relations which are supported only on $a_i$, i.e. live in $S_i$.  These generate $\Lambda_\mathcal{B}'|_{S_i}$, and each is an integer multiple of $\ord(a_i)\mathbf{s}_i$, where $\mathbf{s}_i$ is the $i$-th standard basis vector in $\ZZ^n$.  Write $\alpha_t$ for these integers, $t=1, \ldots, c+1$.

Let us assume heuristically that the values $\alpha_t$ behave as random integers (in the sense of natural density).  The probability that $k$ random integers share no common factor is $1 - 1/\zeta(k)$, where $\zeta$ is the Riemann zeta function.  Therefore the probability that $\Lambda_\mathcal{B}|_{S_i} = \Lambda_\mathcal{B}'|_{S_i}$ would be $1 - 1/\zeta(c+1)$.  

As the convergence $\zeta(n) \rightarrow 1$ as $n \rightarrow \infty$ is exponential, we need only a polynomial number of extra relations to reduce the failure probability to be negligible.
For example, we might expect $\Lambda_\mathcal{B}|_{S_i} = \Lambda_\mathcal{B}'|_{S_i}$ at least 99.9\% of the time if $c \ge 9$.

It is worth noting that it is possible for the algorithm to return $\alpha_t = 0$.  This occurs for elements of $K \subseteq K_i$, which forms a hyperplane of codimension $1$; such occurrences are unlikely (by a cardinality argument similar to that in \cite[Section 2.1]{FonteinWocjan}).

\begin{theorem}
	\label{thm:reduc}
	Let $n$ be an integer, and let $b$ be a function of $n$. 
	Let $\mathcal{B}$ be a factor base of $b$ residues modulo $n$. 
	Suppose an oracle provides random multiplicative relations of size $\le O(\log n)$ amongst the residues $\mathcal{B}$.  
	Then there is a Monte Carlo algorithm to solve the period problem for modulus $n$, whose runtime is $O(b^4 \log b) \poly(\log n)$ and which requires at most $O(b)$ calls to the oracle.  
	Furthermore, under Hypothesis~\ref{heur:sub}, the probability of success approaches $1 - 1/\zeta(c+1)$ where a total of $b+c$ calls are made to the oracle.

\end{theorem}

	Recall that a Monte Carlo algorithm is one whose output may be incorrect with a finite probability.  

\begin{proof}
	The algorithm is as follows (similarly to the latter portions of Algorithm~\ref{alg:factor}).  Let $c > 0$ be an integer.  
	Call the oracle $b+c$ times to obtain $b+c$ relations.
	Isolate variable $a_1$, rewriting the $j$-th relation for each $j$ as
	\[
		a_1^{e_{j,1}} = \prod_{i=2}^b a_i^{e_{j,i}}.
	\]
	Then, by dimensional considerations, there are at least $c+1$ vectors $(c_{k,j})_{j=1}^{b+c}$ for $k=1,\ldots,c+1$ which give linear combinations $\sum_{j=1}^{b+c} c_{k,j} \mathbf{e}_j$ which are supported only on the first entry.  Equivalently, these satisfy
	\[
		a_1^{\sum_{j=1}^{b+c} c_{k,j} e_{j,1}} = 1.
	\]
	We can normalize these relations so the entries $c_{k,j}$ are integers with no common factor, by use of a gcd computation.
	Write $\alpha_k := \sum_{j=1}^{b+c} c_{k,j} e_{j,1}$.
	Let $g$ be the $\gcd$ of the set of exponents $\alpha_k$, $k=1,\ldots,c+1$.  
	Then $g$ is a multiple of the multiplicative order of $a_1$ modulo $n$, say $g = h \ord(a_1)$.

	By the discussion preceding the proof, $h = 
[\Lambda_\mathcal{B}|_{S_i} : \Lambda_\mathcal{B}'|_{S_i}]$.
Under Hypothesis~\ref{heur:sub}, the probability that $h=1$ (and hence the algorithm succeeds) is $1 - 1/\zeta(c+1)$. 

We now take $c$ to be $O(b)$, and consider the runtime.   We must find a basis for the kernel of a $b \times (b+c)$ matrix, where $b+c = O(b)$.  The entries are of size $O(\log n)$ (in the sense that they are integers $\le n$).  Thus the linear algebra takes na\"ive time at most $O(b^4)\poly(\log n)$ resulting in vectors with entries of size $O(b^4)\poly(\log n)$ \cite[Corollary 3.2d]{Schrijver}.  
In the GCD phase, we must perform $O(b)$ gcd operations on integers of size at most $O(b^4\log b)\poly(\log n)$; this has runtime $O(b^4\log b)\poly(\log n)$.
Overall, the runtime of the algorithm is $O(b^4 \log b) \poly(\log n)$.
\end{proof}

\section{Runtime of Algorithm~\ref{alg:factor}}

As usual, we set the notation
\[
	L_x(\alpha,\beta) = \exp( (\beta + o(1)) (\log x)^\alpha (\log\log x)^{1-\alpha}).
\]
We will now show that the algorithm is of runtime $L_n(1/2,\beta)$ for some constant $\beta$, which can be improved by the use of many optimizations developed for the index calculus; see below.  However, since this algorithm is of academic, not practical, interest, we will not devote time to optimizing the constant $\beta$.

The runtime of the linear algebra and gcd phase is given in Theorem~\ref{thm:reduc}.
The relation finding phase is exactly as for the index calculus itself.  If we use the standard notation $u^u$ for the number of trials to find one smooth integer, where $u = \log n/ \log b$, then the runtime is $u^u (b+c) b \pi(b) = u^u O(b^3/\log b)$ with trial division.  Thus, balancing the runtimes we obtain a heuristic runtime of $L_n(1/2,\beta)$ (for the standard analysis of such runtimes, see \cite[Chapter 6]{CrandallPomeranceBook} or \cite[Chapter 15]{GalbraithBook}).

Some relevant implementation notes:

    \begin{enumerate}
	    \item One might include $-1$ in the factor base.
	    \item One may use the linear sieve of Coppersmith, Odlyzko and Schroeppel \cite{CoppersmithOdlzykoSchroeppel}; this improves the relation-finding phase.
	    \item For the linear algebra phase, one might use an algorithm of runtime $O(b^3)\poly(\log n)$ due to Mulders and Storjohann \cite{MuldersStorjohann}.
	    \item One might use the elliptic curve method \cite{Lenstra} to remove all prime factors below a bound before attempting this algorithm.  Having many small factors will result in the order of $g$ being smaller more often, which gives a higher failure rate in using the order to obtain the factorization.
	    \item We might test for the existence of a non-trivial kernel periodically as we generate relations, since if the order of $g$ is small, the algorithm actually requires fewer relations.
	    \item In the same vein, we might use a single kernel element, when it is found, to obtain a multiple of $\ord(g)$ and then do further linear algebra modulo that modulus.
	    \item We may need to return to the relation-finding and add more relations if we do not find enough elements of the kernel which result in $\alpha_k \neq 0$.
	    \item Many implementation tricks for the index calculus may apply in this situation, see for example \cite{Odlyzko84}.
	    \item The number field sieve as in \cite[Section 3.1]{Gordon} can be adapted to speed up the relation-finding phase.
    \end{enumerate}

    \section{Example}

    Let us compute the order of $g = 43$ modulo $n = 62389$ and use this to factor $n$.  We will use $B = 50$, resulting in a factor base of $b=15$ primes $2 \le p \le 47$.  We will need $25$ relations (this assumes $c=10$).  With $188$ smoothness tests, we find the relations:

    \begin{minipage}{0.29\textwidth}
    \begin{align*}
43 ^{ 55571 }&= 2^{3} \cdot 3^{3} \cdot 7 \cdot 29 ,\\
43 ^{ 51344 }&= 5^{4} ,\\
43 ^{ 1724 }&= 2 \cdot 5^{3} \cdot 7 \cdot 23 ,\\
43 ^{ 9399 }&= 3 \cdot 13 \cdot 37 ,\\
43 ^{ 56136 }&= 2 \cdot 3 \cdot 11^{2} \cdot 13 ,\\
43 ^{ 53393 }&= 5^{4} \cdot 41 ,\\
43 ^{ 24567 }&= 2^{4} \cdot 7 \cdot 23^{2} ,\\
43 ^{ 2484 }&= 2 \cdot 3^{2} \cdot 13 \cdot 37 ,\\
43 ^{ 39818 }&= 7^{2} ,\\
	    \end{align*}
    \end{minipage}\hfill
    \begin{minipage}{0.29\textwidth}
	    \begin{align*}
43 ^{ 41451 }&= 2^{2} \cdot 5 \cdot 7 \cdot 11^{2} ,\\
43 ^{ 53596 }&= 3^{3} \cdot 11 \cdot 43 ,\\
43 ^{ 12688 }&= 2^{3} \cdot 3 \cdot 7 \cdot 19^{2} ,\\
43 ^{ 10480 }&= 2^{3} \cdot 3^{3} \cdot 5 \cdot 13 ,\\
43 ^{ 19831 }&= 2^{8} \cdot 3 \cdot 5 \cdot 11 ,\\
43 ^{ 27853 }&= 2^{6} \cdot 3^{2} \cdot 5 \cdot 7 ,\\
43 ^{ 25154 }&= 2^{5} \cdot 31 \cdot 37 ,\\
43 ^{ 9481 }&= 2^{3} \cdot 7 \cdot 11 ,\\
\end{align*}
    \end{minipage}\hfill
    \begin{minipage}{0.29\textwidth}
	    \begin{align*}
43 ^{ 20 }&= 2^{2} \cdot 5^{3} \cdot 7^{2} ,\\
43 ^{ 25418 }&= 2^{5} \cdot 3 \cdot 17 \cdot 19 ,\\
43 ^{ 50821 }&= 5^{2} \cdot 41 ,\\
43 ^{ 46106 }&= 2 \cdot 3 \cdot 7 \cdot 11^{2} ,\\
43 ^{ 14141 }&= 2 \cdot 3 \cdot 5^{2} \cdot 7 \cdot 19 ,\\
43 ^{ 26246 }&= 2 \cdot 3^{3} \cdot 5 \cdot 41 ,\\
43 ^{ 10795 }&= 2 \cdot 5^{3} \cdot 7 \cdot 11 ,\\
43 ^{ 20889 }&= 5 \cdot 11 \cdot 37 ,\\
\end{align*}
    \end{minipage}

The relation matrix is (cols are relations):
\[ \scriptsize
\left(\begin{array}{rrrrrrrrrrrrrrrrrrrrrrrrr}
3 & 0 & 1 & 0 & 1 & 0 & 4 & 1 & 0 & 2 & 0 & 3 & 3 & 8 & 6 & 5 & 3 & 2 & 5 & 0 & 1 & 1 & 1 & 1 & 0 \\
3 & 0 & 0 & 1 & 1 & 0 & 0 & 2 & 0 & 0 & 3 & 1 & 3 & 1 & 2 & 0 & 0 & 0 & 1 & 0 & 1 & 1 & 3 & 0 & 0 \\
0 & 4 & 3 & 0 & 0 & 4 & 0 & 0 & 0 & 1 & 0 & 0 & 1 & 1 & 1 & 0 & 0 & 3 & 0 & 2 & 0 & 2 & 1 & 3 & 1 \\
1 & 0 & 1 & 0 & 0 & 0 & 1 & 0 & 2 & 1 & 0 & 1 & 0 & 0 & 1 & 0 & 1 & 2 & 0 & 0 & 1 & 1 & 0 & 1 & 0 \\
0 & 0 & 0 & 0 & 2 & 0 & 0 & 0 & 0 & 2 & 1 & 0 & 0 & 1 & 0 & 0 & 1 & 0 & 0 & 0 & 2 & 0 & 0 & 1 & 1 \\
0 & 0 & 0 & 1 & 1 & 0 & 0 & 1 & 0 & 0 & 0 & 0 & 1 & 0 & 0 & 0 & 0 & 0 & 0 & 0 & 0 & 0 & 0 & 0 & 0 \\
0 & 0 & 0 & 0 & 0 & 0 & 0 & 0 & 0 & 0 & 0 & 0 & 0 & 0 & 0 & 0 & 0 & 0 & 1 & 0 & 0 & 0 & 0 & 0 & 0 \\
0 & 0 & 0 & 0 & 0 & 0 & 0 & 0 & 0 & 0 & 0 & 2 & 0 & 0 & 0 & 0 & 0 & 0 & 1 & 0 & 0 & 1 & 0 & 0 & 0 \\
0 & 0 & 1 & 0 & 0 & 0 & 2 & 0 & 0 & 0 & 0 & 0 & 0 & 0 & 0 & 0 & 0 & 0 & 0 & 0 & 0 & 0 & 0 & 0 & 0 \\
1 & 0 & 0 & 0 & 0 & 0 & 0 & 0 & 0 & 0 & 0 & 0 & 0 & 0 & 0 & 0 & 0 & 0 & 0 & 0 & 0 & 0 & 0 & 0 & 0 \\
0 & 0 & 0 & 0 & 0 & 0 & 0 & 0 & 0 & 0 & 0 & 0 & 0 & 0 & 0 & 1 & 0 & 0 & 0 & 0 & 0 & 0 & 0 & 0 & 0 \\
0 & 0 & 0 & 1 & 0 & 0 & 0 & 1 & 0 & 0 & 0 & 0 & 0 & 0 & 0 & 1 & 0 & 0 & 0 & 0 & 0 & 0 & 0 & 0 & 1 \\
0 & 0 & 0 & 0 & 0 & 1 & 0 & 0 & 0 & 0 & 0 & 0 & 0 & 0 & 0 & 0 & 0 & 0 & 0 & 1 & 0 & 0 & 1 & 0 & 0 \\
0 & 0 & 0 & 0 & 0 & 0 & 0 & 0 & 0 & 0 & 1 & 0 & 0 & 0 & 0 & 0 & 0 & 0 & 0 & 0 & 0 & 0 & 0 & 0 & 0 \\
0 & 0 & 0 & 0 & 0 & 0 & 0 & 0 & 0 & 0 & 0 & 0 & 0 & 0 & 0 & 0 & 0 & 0 & 0 & 0 & 0 & 0 & 0 & 0 & 0
\end{array}\right)
\]
The following matrix is made up of rows representing the right kernel:
\[ \scriptsize
\left(\begin{array}{rrrrrrrrrrrrrrrrrrrrrrrrr}
0 & 1 & 0 & 0 & 0 & 0 & 0 & 0 & 0 & 0 & 0 & 0 & 0 & 0 & 12 & 0 & -23 & 3 & 0 & 14 & 18 & 0 & -14 & -13 & 0 \\
0 & 0 & 2 & 0 & 0 & 0 & -1 & 0 & 0 & 0 & 0 & 0 & 0 & 0 & 7 & 0 & -12 & 1 & 0 & 8 & 10 & 0 & -8 & -8 & 0 \\
0 & 0 & 0 & 1 & 0 & 0 & 0 & 0 & 0 & 0 & 0 & 0 & -1 & 0 & 0 & 0 & 1 & -1 & 0 & -1 & -1 & 0 & 1 & 2 & -1 \\
0 & 0 & 0 & 0 & 1 & 0 & 0 & 0 & 0 & 0 & 0 & 0 & -1 & 0 & 8 & 0 & -15 & 2 & 0 & 8 & 10 & 0 & -8 & -7 & 0 \\
0 & 0 & 0 & 0 & 0 & 1 & 0 & 0 & 0 & 0 & 0 & 0 & 0 & 0 & 13 & 0 & -25 & 3 & 0 & 14 & 19 & 0 & -15 & -13 & 0 \\
0 & 0 & 0 & 0 & 0 & 0 & 0 & 1 & 0 & 0 & 0 & 0 & -1 & 0 & 4 & 0 & -7 & 0 & 0 & 4 & 5 & 0 & -4 & -2 & -1 \\
0 & 0 & 0 & 0 & 0 & 0 & 0 & 0 & 1 & 0 & 0 & 0 & 0 & 0 & 3 & 0 & -5 & -1 & 0 & 3 & 3 & 0 & -3 & -1 & 0 \\
0 & 0 & 0 & 0 & 0 & 0 & 0 & 0 & 0 & 1 & 0 & 0 & 0 & 0 & 8 & 0 & -16 & 2 & 0 & 9 & 11 & 0 & -9 & -8 & 0 \\
0 & 0 & 0 & 0 & 0 & 0 & 0 & 0 & 0 & 0 & 0 & 1 & 0 & 0 & 2 & 0 & -5 & 1 & 0 & 2 & 3 & -2 & -2 & -1 & 0 \\
0 & 0 & 0 & 0 & 0 & 0 & 0 & 0 & 0 & 0 & 0 & 0 & 0 & 1 & 6 & 0 & -15 & 3 & 0 & 8 & 11 & 0 & -8 & -8 & 0 \\
0 & 0 & 0 & 0 & 0 & 0 & 0 & 0 & 0 & 0 & 0 & 0 & 0 & 0 & 14 & 0 & -27 & 3 & 0 & 16 & 20 & 0 & -16 & -13 & 0
\end{array}\right)
\]
The corresponding $\alpha_k$ are:
\[
1201200, \; 631400,  \;-61600,  \;708400,  \;1232000,  \;323400,  \;277200,  \;754600,  \;169400,  \;662200,  \;1309000.
\]
Their gcd is $15400$.
 We check that
\[
	43^{15400} = 1, \quad 43^{15400/2} = 51174 \neq  \pm 1 \pmod{n}.
\]
and therefore taking
\[
	\gcd(51174-1, 62389) = 701
\]
reveals a non-trivial factor.  In fact, $ 62389= 701 \cdot 89$.

\bibliographystyle{abbrv}  
\bibliography{index-factor} 

\end{document}